\documentclass[12pt,reqno]{amsart}
\usepackage{mathrsfs,cmtiup}
\usepackage{hyperref}
\usepackage{amssymb,amsthm}
\usepackage{amsfonts,cmtiup,comment,stmaryrd}

\makeatletter
\def\blfootnote{\xdef\@thefnmark{}\@footnotetext}
\makeatother

\newtheorem{theorem}{Theorem}[section]
\newtheorem{lemma}[theorem]{Lemma}
\newtheorem{proposition}[theorem]{Proposition}

\newcommand{\f}{\varphi}

\newcommand{\Z}{\Bbb Z}

\newcommand{\w}{\omega }

\newcommand{\ka}{\varkappa }
\renewcommand{\leq}{\leqslant}
\renewcommand{\geq}{\geqslant}
\theoremstyle{definition}

\newtheorem*{definition*}{Definition}
\numberwithin{equation}{section}

\begin{document}

\title{Profinite groups with coprime automorphisms\\ whose orbits generate nilpotent subgroups}

\author{Evgeny Khukhro}
\address{E. I. Khukhro: Charlotte Scott Research Centre for Algebra, University of Lincoln, U.K.}
\email{khukhro@yahoo.co.uk}

\author{Pavel Shumyatsky}

\address{P. Shumyatsky: Department of Mathematics, University of Brasilia, DF~70910-900, Brazil}
\email{pavel@unb.br}

\thanks{The first author was partially supported by the International Center for Mathematics at SUSTech in Shenzhen. The second author was supported by  FAPDF and CNPq.}


\begin{abstract}  Suppose that a profinite group $G$ admits a finite group of coprime automorphisms $A$ such that for any $\f\in A$ every subgroup $\langle x^{\langle\f\rangle}\rangle$ generated by the $\langle\f\rangle$-orbit of an element $x\in G$ is nilpotent. It is proved that $[G,A]$ is pronilpotent, and if $A$ is soluble and the fixed-point subgroup $C_G(A)$ is locally nilpotent, then $G$ is  locally nilpotent.
\end{abstract}

\maketitle

\section{Introduction}
A continuous automorphism  $\varphi$  of a profinite group $G$ is called a \textit{coprime automorphism} of $G$ if $\f$ has finite order and $G$ is an inverse limit of finite groups of order coprime to the order of $\varphi$.
Many properties of coprime automorphisms of profinite groups are similar to the properties of coprime automorphisms of finite groups. For example, if $\varphi$ is a coprime automorphism of $G$, then for any (closed) normal $\varphi$-invariant subgroup $N$ the fixed points of the induced automorphism (denoted by the same letter) in $G/N$ are images of the fixed points in $G$, that is, $C_{G/N}(\varphi )=C_G(\varphi )N/N$. Therefore, if $\varphi$ is a  coprime automorphism  of prime order $q$ such that $C_G(\varphi )=1$, Thompson's theorem  \cite{tho} implies that $G$ is pronilpotent. Moreover, $G$ is nilpotent and the nilpotency class is bounded in terms of $q$ by Higman's theorem~\cite{hig}.
These results can be viewed as a starting point of an important direction in the study of profinite groups where conditions on automorphisms and their fixed points are used for obtaining various nilpotency-type results; see, for example, \cite{aks, acc-shu,kh-sh22, shu00a}.

In this paper we study profinite groups $G$ admitting a finite group of  coprime automorphism $A$ such that the orbits of elements of $G$ under cyclic subgroups of $A$ generate nilpotent subgroups.

\begin{theorem}\label{t-2}
Suppose that a profinite group $G$ admits a finite group $A$ of coprime automorphisms such that for any $\f\in A$  every subgroup $\langle x^{\langle\f\rangle}\rangle$ generated by the $\langle\f\rangle$-orbit of an element $x\in G$ is nilpotent.
\begin{itemize}
  \item[{\rm (a)}] Then the commutator subgroup $[G,A]$ is pronilpotent.
  \item[{\rm (b)}] If in addition $A$ is soluble and the fixed-point subgroup $C_G(A)$ is locally nilpotent, then $G$ is  locally nilpotent.
\end{itemize}
\end{theorem}

Examples show that in part (a) one cannot claim that the whole group $G$ is pronilpotent, and in part (b) one cannot drop either the local nilpotency condition on $C_G(A)$ or the nilpotency condition on the subgroups  $\langle x^{\langle\f\rangle}\rangle$.

Firstly Theorem~\ref{t-2} is proved for the case where $G$ is finite (Theorem~\ref{t-fin}), which is equivalent to proving part (a) and the pronilpotency of $G$ in part (b). The proof of Theorem~\ref{t-fin} uses the classification of finite simple groups. Once $G$ is pronilpotent, it remains to prove that if $G$ is finitely generated, then all Sylow $p$-subgroups $G_p$ of $G$ are nilpotent with a uniform upper bound for the nilpotency class, independent of $p$. Initially we prove that every $G_p$ is nilpotent, with a bound for the nilpotency class depending on $p$. Then a bound independent of $p $ is obtained for all sufficiently large primes $p$. In both cases we apply Lie ring methods, and the crucial tool is Zelmanov's theorem \cite{ze92,ze95,ze17} on Lie algebras and some of its consequences.
Other important ingredients include  criteria for a pro-$p $ group to be $p $-adic analytic in terms of the associated Lie algebra and in terms of bounds for the rank due to Lazard~\cite{laz} and Lubotzky and Mann \cite{lu-ma}, and a theorem of Bahturin and Zaicev~\cite{ba-za} on Lie algebras admitting a group of automorphisms whose fixed-point subalgebra satisfies a polynomial identity.

\section{Preliminaries}\label{s-prel}

We shall routinely say that a group or a Lie ring is nilpotent  of class $k$ if it is nilpotent of class at most $k$. We write ``$(a,b\dots)$-bounded'' to abbreviate ``bounded above by some function depending only on the parameters $a,b\dots$".

\subsection*{\bf Lie rings and algebras}
We call products in Lie rings and algebras commutators. Simple commutator notation is used for left-normed commutators  $[x_1,\dots  , x_k] = [...[x_1, x_2],\dots ,  x_k]$, and Engel commutators are denoted by $[x,{}_ny]=[x,y,y,\dots, y]$, where $y$ occurs $n$ times. An element $a$ of a Lie ring or a Lie algebra $L$ is said to be {ad-nilpotent} of index $n$  if  $[x,{}_n a] = 0$ for all $x\in L$; here  $n$ does not have to be the least integer with this property, so $a$ is also ad-nilpotent of index $m$ for any $m\geq n$.

The next theorem is a deep result of Zelmanov \cite{ze92,ze95,ze17}.

\begin{theorem}\label{t-z}
 Let $L$ be a Lie algebra over a field and suppose that $L$ satisfies
a polynomial identity. If $L$ can be generated by a finite set $X$ such that every
commutator in elements of $X$ is ad-nilpotent, then $L$ is nilpotent.
\end{theorem}

Henceforth,  elements of $X$ are also considered to be commutators in elements of $X$.

An important criterion for a Lie algebra to satisfy a polynomial identity
is provided by the next theorem, which was proved by Bahturin and Zaicev \cite{ba-za}  for
soluble group of automorphisms and later extended by Linchenko \cite{lin} to the general case. We use the centralizer notation for the fixed point subalgebra $C_L(A)$ of a group of automorphisms $A$ of~$L$.

\begin{theorem}\label{t-bz}
Let $L$ be a Lie algebra over a field $K$. Assume that a finite
group $A$ acts on $L$ by automorphisms in such a manner that $C_L(A)$ satisfies
a polynomial identity. Assume further that the characteristic of $K$ is either $0$
or coprime with the order of $A$. Then $L$ satisfies a polynomial identity.
\end{theorem}

Both Theorems~\ref{t-z} and \ref{t-bz} admit respective quantitative versions (see for
example \cite{shu00}). For our purposes, we shall need the following result for Lie rings proved in \cite[Proposition~2.6]{danilo}, which combines both versions. As usual, $\gamma _i(L)$ denotes the $i$-th
term of the lower central series of $L$.

\begin{proposition}\label{p-danilo}
Let $L$ be a Lie ring and $A$ a finite group of automorphisms of
$L$ such that $C_L(A)$ satisfies a polynomial identity $f\equiv 0$. Suppose that
$L$ is generated by an $A$-invariant set of $m$ elements such that every commutator in these elements is ad-nilpotent of index~$n$. Then there exist positive
integers $e$ and $c$ depending only on $|A|$, $f$,  $m$, and $n$ such that $e\gamma _c(L) = 0$.
\end{proposition}

The following folklore lemma simplifies verification of ad-nilpotency of an element in a finitely generated  Lie algebra over a field of positive characteristic.

\begin{lemma}
  \label{l-charp}
  Suppose that $L$ is a Lie algebra over a field of characteristic $p\ne 0$ finitely generated by elements $g_1,\dots ,g_n$. If for an element $h\in L$ there are positive integers $k_1,\dots,k_n$ such that $[g_i,{}_{k_i}h]=0$ for all $i=1,\dots,n$, then $h$ is an ad-nilpotent element of $L$ of index bounded in terms of $p$ and $\max\{k_1,\dots ,k_n\}$.
\end{lemma}

\begin{proof}
Choose a positive integer $s$ such that  $p ^s\geqslant \max \{k_1,\dots ,k_m\}$. Then $[g_j,{}_{p ^s}h]=0$ for all $j$. In characteristic $p $ this implies that
\begin{equation}\label{e-qs}
[\varkappa ,{}_{p ^s}h]=0
\end{equation}
for any commutator $\varkappa$ in the $g_i$. This follows by induction on the weight of $\varkappa$ from the formula
$$
[[u,v],{}_{p ^s}w]=[[u,{}_{p ^s}w],v]+[u,[v,{}_{p ^s}w]],
$$
which holds in any Lie algebra of characteristic $p$. This formula is a consequence of the Leibnitz rule
$$
[[u,v],{}_{n}w]=\sum_{i=0}^n {n\choose i}
[[u,{}_{i}w],\,[v,{}_{n-i}w]]
$$
(where $[a,{}_{0}b]=a$), because the binomial coefficient ${p ^s\choose i}$ is divisible by $p $ unless $i=0$ or $i=p ^s$.
Since any element of  $L$ is a linear combination of commutators in the $g_i$, equation \eqref{e-qs} by linearity implies that $h$ is ad-nilpotent of index~$p ^s$.
\end{proof}

The following lemma appeared in \cite[Lemma~5]{kh-sh99}, where it was stated for Lie algebras, so we repeat the proof for Lie rings for the benefit of the reader.

\begin{lemma}\label{l-be}
Let $L$ be a Lie ring, and $M$ a nilpotent subring of class $c$ generated by finitely many elements $g_1,\dots,g_m$ such that all commutators in these elements are ad-nilpotent of index $n$ in $L$. Then for some number $\varepsilon =\varepsilon (m,n,c)$ bounded in terms of $m$, $n$, $c$ we have $[L,\,\underbrace{M, M,\dots ,M}_\varepsilon ]=0$. In particular, every element of $M$ is ad-nilpotent of index $\varepsilon (m,n,c)$.
\end{lemma}

\begin{proof}
Since $M=\langle g_1,\dots , g_m\rangle$ is nilpotent of class $c$,
there are only $(c,m)$-boundedly many $N=N(c,m)$ nonzero commutators in the $g_i$. We enumerate all these commutators as $\ka _1,\ka_2,\dots,\ka_N$ ordered in accordance with the increase of their weight (and arbitrarily among commutators of the same weight). It is sufficient to show that for some $(m,n,c)$-bounded number $\varepsilon =\varepsilon (m,n,c)$ for any $u\in  L$ and any $v_i\in M$ we have
  $[u,v_1,\dots,v_\varepsilon]=0$. We express all the $v_i$ as Lie ring words in the generators $g_i$ and repeatedly use the linearity and the Jacobi identity $[[a,b],c]=[[a,c],b]+[a,[b,c]]$ in a collection process to  represent $[u,v_1,\dots,v_\varepsilon]$ as a linear combination of commutators of the form
\begin{equation}\label{e-long}
[[[u,{}_{k_1}\ka _1],{}_{k_2}\ka _2],\,\dots, {}_{k_N}\ka _N].
\end{equation}
 The weights of the $\ka_i$  as commutators in the $g_i$ are at most $c$, while the total sum of these weights with multiplicities  is at least $\varepsilon$. Since $N$ is $(c,m)$-bounded, there is a $(m,n,c)$-bounded  $\varepsilon $ such that at least one of the $k_i$ must be at least $n$. Then  the commutator \eqref{e-long} is equal to 0, since $\ka _i$ is ad-nilpotent of index $n$. \end{proof}

\subsection*{\bf Profinite groups and their automorphisms}
A profinite group is always regarded as a topological group. Unless stated otherwise, a~subgroup of a topological group will always mean a closed subgroup, all homomorphisms are continuous, and quotients are by closed normal subgroups. This also applies to taking commutator subgroups, normal closures, subgroups generated by subsets, etc. We say that a subgroup is generated by a subset $X$ if it is generated by $X$ as a topological group. Note that if $\varphi$ is a continuous automorphism of a topological group $G$, then the fixed-point subgroup $C_G(\varphi )$  is closed.

Recall that a {pronilpotent} group is a pro-(finite nilpotent) group, that is, an inverse limit of finite nilpotent groups. For a prime $p $, a pro-$p $ group is  an inverse limit of finite $p $-groups. The Frattini subgroup $P'P^p $  of a finitely generated pro-$p $ group $P$ is the product of the derived subgroup $P'$ and the subgroup generated by all $p$-th powers of elements of $P$; it is also known as the dimension subgroup $D_2(P)$. A subset  generates  $P$ (as a topological group) if and only if its image  generates the elementary abelian quotient $P/(P'P^p )$. See, for example, \cite{wil} for these and other properties of profinite groups.

 The following lemma collects several well-known properties of coprime automorphisms of finite and profinite groups. For finite groups see, for example, \cite[Ch.~5,~6]{gor}, and the profinite versions are straightforward using the inverse limit argument; see for example, \cite[Proposition~2.3.16]{rib-zal}. This lemma will often be used without explicit references.

\begin{lemma}\label{l-nakr}
Let $A$ be a group of coprime automorphisms of a profinite group $G$.
\begin{itemize}
  \item[\rm (a)]  If $N$ is an $A$-invariant closed normal  subgroup of~$G$, then the fixed points of $A$ in  the quotient $G/N$ are covered by the fixed points of $A$ in $G$, that is, $C_{G/N}(A) = C_G(A)N/N$.

  \item[\rm (b)] We have $[[G, A],A]=[G,A]$.

  \item[\rm (c)]   If $G$ is abelian, then $G=[G,A]\times C_G(A)$.

  \item[\rm  (d)] If $A$ is abelian noncyclic, then $G=\prod_{1\ne a\in A}C_G(a)$.

  \item[\em (e)] For every prime $p$ there is an $A$-invariant Sylow $p$-subgroup of $G$.
\end{itemize}
  \end{lemma}

 Suppose that $A$ is a group of coprime automorphisms of a profinite group $G$. Then $A$ induces a group of coprime automorphisms on every quotient of $G$ by a closed normal $A$-invariant subgroup. Suppose that the fixed-point subgroup $C_G(A)$ is locally nilpotent. Then by Lemma~\ref{l-nakr} this property is inherited by every quotient of $G$ by an $A$-invariant normal subgroup for the induced group of automorphisms.

The following useful lemma was proved for coprime automorphisms of finite nilpotent groups in \cite[Lemma~2.4]{rod-shu}; for non-nilpotent groups the lemma fails (as shown by an example in \cite{rod-shu}).

\begin{lemma}\label{l-copr}
Let $\varphi$ be a coprime automorphism of a finite nilpotent
group $G$. Then any element $g\in G$ can be uniquely written in the form
$g=c[u,\varphi ]$, where $c\in C_G(\varphi)$ and $u\in G$.
 \end{lemma}

A similar result follows also for pronilpotent groups.

\begin{lemma}\label{l-coprp}
Let $\varphi$ be a coprime automorphism of a pronilpotent
group $G$. Then any element $g\in G$ can be uniquely written in the form
$g=c[u,\varphi ]$, where $c\in C_G(\varphi)$ and $u\in G$.
 \end{lemma}

 \begin{proof}
  The set $I_G(\f)=\{[x,\f ]\mid x\in G\}$ is closed as the image of the continuous mapping $x\mapsto  [x,\f ]$. The set $C_G(\varphi)$ is also closed. Hence the product
  $$
  C_G(\varphi)I_G(\f)=\{c[x,\f]\mid c\in C_G(\f),\;x\in G\}
   $$
   is also closed, since $G$ is compact. Every normal open subgroup $K$ contains a $\f$-invariant normal open subgroup $\bigcap_{i=0}^{|\f|-1}K^{\f^i}$; therefore $\f$-invariant normal open subgroups form a basis of neighbourhoods of the identity in $G$. In every finite quotient by a $\f$-invariant closed normal subgroup the image of the closed set $C_G(\varphi)I_G(\f)$ is equal to the whole quotient by Lemma~\ref{l-copr}.  Hence this product is equal to $G$. We now prove uniqueness by contradiction: suppose that $c_1[x_1,\f]=c_2[x_2,\f]$ for $c_1,c_2\in C_G(\f)$ such that $c_1\ne c_2$. Then there is a $\f$-invariant normal open subgroup $K$ such that the images of $c_1$ and $c_2$ in $G/K$ are different. This contradicts the uniqueness in Lemma~\ref{l-copr}.
 \end{proof}
\begin{lemma}

\label{l-cleft1}
Suppose that $\f$ is a coprime automorphism  of a pronilpotent  group $G$.  Let $g\in G$  and $c\in C_G(\f)$. Then the subgroup $\langle (c[g,\f])^{\langle\f\rangle}\rangle$ generated by the $\langle\f\rangle$-orbit of $c[g,\f]$ coincides with the subgroup $\langle c,[g,\f]^{\langle\f\rangle}\rangle$ generated by $c$ and the $\langle\f\rangle$-orbit of $[g,\f]$. \end{lemma}

\begin{proof}
Let  $K=\langle (c[g,\f])^{\langle \f\rangle}\rangle$.  By Lemma~\ref{l-coprp} there are elements $v\in C_K(\f)$ and  $u\in K$ such that $c[g,\f]=v[u,\f]$. By  Lemma~\ref{l-coprp} such a representation is  unique, so that $c=v\in K$ and $[g,\f]=[u,\f ]\in K$. Thus both $c$ and $[g,\f]$  are contained in the $\f$-invariant subgroup $\langle (c[g,\f])^{\langle\f\rangle}\rangle$. Therefore, $\langle c,[g,\f]^{\langle\f\rangle}\rangle\leq \langle (c[g,\f])^{\langle \f\rangle}\rangle$. The reverse inclusion is obvious.
  \end{proof}

Let $\f$ be  a coprime automorphism of a profinite group $G$.
Suppose that every subgroup $\langle x^{\langle\f\rangle}\rangle$ generated by the $\langle\f\rangle$-orbit of an element $x\in G$ is nilpotent. Then this property is inherited by every quotient of $G$ by a $\f$-invariant normal subgroup for the induced automorphism, as well as by any power $\f^k$ of the automorphism $\f$.

\section{Pronilpotency}

In this section we prove part (a) of Theorem~\ref{t-2}, that $[G,A]$ is pronilpotent, as well as the pronilpotency of $G$ in part (b), when $A$ is soluble and $C_G(A)$ is locally nilpotent. Both these results are essentially about finite groups.

\begin{theorem}\label{t-fin}
Suppose that a finite group $G$ admits a group $A$ of coprime automorphisms such that for every $\f\in A$ every subgroup $\langle x^{\langle\f\rangle}\rangle$ generated by the $\langle\f\rangle$-orbit of an element $x\in G$ is nilpotent.
\begin{itemize}
  \item[\rm (a)] Then the commutator subgroup $[G,A]$ is nilpotent.
  \item[\rm (b)] If in addition $A$ is soluble and the fixed-point subgroup $C_G(A)$ is nilpotent, then $G$ is  nilpotent.
\end{itemize} \end{theorem}

\begin{proof}
(a) Since $[G,A]=\prod_{\f\in A}[G,\f]$, it is sufficient to prove that each of the normal subgroups $[G,\f]$ is nilpotent. The hypotheses are satisfied by the group of automorphisms~$\langle\f\rangle$ for any $\f\in A$. Therefore we can assume that $A=\langle\f\rangle$ is cyclic. Arguing by contradiction, let $G$ be a counterexample of minimal order. Then $G=[G,\f]$. Let $M$ be a minimal $\f$-invariant normal subgroup of $G$. First suppose that $M=G$. Then $M$ is non-abelian, so that $G=S_1\times \cdots \times S_k$ is a direct product of isomorphic  non-abelian simple groups transitively permuted by $\f$. If $k=1$, then $G$ is simple and there exists $1\ne x\in G$ such that
 $\langle x,\f\rangle=\langle x^{\langle\f\rangle}\rangle\rtimes \langle \f\rangle=G\rtimes \langle \f\rangle$ by the ``3/2 Theorem" of Guralnick and Kantor \cite{gur-kan}, contrary to $\langle x^{\langle\f\rangle}\rangle$ being nilpotent by hypothesis. Now let $k\geq 2$; we can assume without loss of generality that $S_1^\f=S_2$. Then for $g=(a_1,a_2,\dots )\in G$ with $a_i\in S_i$, we have $g^\f=(a_1,a_2,\dots )^\f=(*,a_1^\f,\dots )$. The subgroup $\langle g,g^\f\rangle$ is nilpotent by hypothesis, but the second `coordinates' $a_2$ and $a_1^\f$ of these elements  can be chosen to generate $S_2$, since every finite simple group is 2-generated, as follows from the classification of finite simple groups; a contradiction.

 Thus, we can assume that $M<G$. By the minimality, both $G/M$ and $[M,\f]$ are nilpotent. If $\f$ centralizes $M$, then $G=[G,\f]$ also centralizes $M$ and then $G$ is nilpotent. Hence,  $[M,\f]\ne 1$. If $M$ is non-abelian, then $[M,\f]$ is a direct product of non-abelian simple groups, a contradiction. Hence $M$ is an elementary abelian $p$-group for some prime $p$.

 Since $G$ is not nilpotent by assumption, there is a Sylow $q$-subgroup $Q$ for $q\ne p$ that acts nontrivially on $M$; we can choose $Q$ to be $\f$-invariant by Lemma~\ref{l-nakr}(e), so that $Q=[Q,\f]$ (since $G=[G,\f]$). Note that $Q$ does not centralize $[M,\f]$: otherwise $M=[M,\f]\times N$ for a $Q\langle\f\rangle$-invariant subgroup $N$ contained in $C_M(\f)$ and then $Q=[Q,\f]$ also centralizes $N$ and therefore acts trivially on $M$ contrary to the choice of $Q$. It follows that the  $\f$-invariant subgroup $[MQ,\f]$ is also non-nilpotent, and thus $G=MQ$ by the minimality. Since $MQ/C_Q(M)$ and $MQ/C_M(Q)$ are also non-nilpotent, we have $C_Q(M)=1$ and $C_M(Q)=$ by the minimality.

 Let $Q_1$ be the minimal $\f$-invariant subgroup of $Q$ on which $\f$ acts nontrivially; this is a special $q$-group and $Q_1=[Q_1,\f]$. Since $Q_1$ acts faithfully on $M$, the above arguments applied to $Q_1$ show that $Q=Q_1$ by the minimality.

 If $Q$ is non-abelian, then $C_Q(\f )=Z(Q)\ne 1$. For any  $c\in Z(Q)$ and any $m\in M$ the subgroup $\langle cm, cm^\f,\dots \rangle$ is nilpotent. In particular, this subgroup contains $m^{-1}m^\f=[m,\f]$ and $cm$, so that
 $$
 1=[[m,\f], cm,cm,\dots ]=[[m,\f], c,c,\dots ],
  $$
  whence $[[M,\f], Z(Q)]=1$, since the orders of $Z(Q)$ and $M$ are coprime. Since $C_M(Z(Q))$ is a nontrivial $\f$-invariant normal subgroup of $G$, we have   $C_M(Z(Q))=M$ by the minimality, a contradiction, since $Q$ acts faithfully on $M$.  Hence   $Q$ is elementary abelian, $C_Q(\f)=1$, and $Q=\langle g^{\langle\f\rangle}\rangle$ for any $1\ne  g\in Q$.

 For any $m\in M$ and $1\ne  g\in Q$, the subgroup $H=\langle gm, g^\f m^\f,\dots \rangle$ is nilpotent. The image $HM/M$ is equal to $QM/M$, and therefore $H$ contains a Sylow $q$-subgroup of $G$. In fact, $H$ must be equal to a Sylow $q$-subgroup, since otherwise $H \cap M$ is a nontrivial $p$-subgroup centralizing a Sylow $q$-subgroup of $G$ and hence $C_M(Q)\ne 1$, a contradiction. Thus $H$ is a $q$-group. It follows that $|mg|=q$ for any $m\in M$ and $1\ne  g\in Q$. Then $Q$ must be cyclic (of order $q$), since otherwise $M=\prod_{1\ne g\in Q}C_M(g)$ by Lemma~\ref{l-nakr}(d) and there would exist  $m\in M$ and $1\ne  g\in Q$ with $|mg|=pq$.

 Let $Q=\langle g\rangle$. For any $m\in M$ the subgroup $\langle g^m, (g^\f)^{m^\f},\dots \rangle$ is nilpotent. Hence both elements $g^m$, $(g^\f)^{m^\f}$ of order $q$ are in the same cyclic subgroup of order $q$, that is,
 $\langle g^m\rangle=\langle (g^\f)^{m^\f}\rangle=\langle g^{m^\f}\rangle$, since $Q$ is $\f$-invariant. Then $m^{-1}m^\f$ normalizes and therefore centralizes~$Q$. Thus, $[M,\f]$ commutes with $Q$, a contradiction.

 (b)   First we note that it is sufficient to prove the theorem for the case where $A$ is cyclic. Indeed,  since $A$ is soluble, there is a normal subgroup $A_1$ of prime index. Let $\f\in A\setminus A_1$. Then  $C_G(A_1)$ admits the cyclic group of automorphisms induced by $\langle\f\rangle$, which satisfies on $C_G(A_1)$ the same hypotheses as $A$. Assuming that the theorem is proved in the case of cyclic $A$, we obtain that $C_G(A_1)$ is nilpotent. Now the group $G$ admits the group of automorphisms $A_1$ of smaller order satisfying the same hypotheses, and therefore the proof concludes by the obvious induction on $|A|$.   Thus, from now on we assume that $G=\langle\f\rangle$ is cyclic.

 We argue by contradiction and let $G$ be a counterexample of minimal order. The subgroup $[G,\f]$ is nilpotent by part (a), and $G/[G,\f]=C_G(\f)[G,\f]/[G,\f]$ is nilpotent by hypothesis. Since $G$ is not nilpotent, there is a Sylow $q$-subgroup $Q$ of $C_G(\f)$ that acts nontrivially on a Sylow $p$-subgroup $P$ of $[G,\f]$ for $p\ne q$. Note that both $Q$ and $P$ are $\f$-invariant. Hence $G=P Q$ by the minimality. Since $Q$ acts also nontrivially on the Frattini quotient of $P$, we must have $\Phi (P)=1$ by the minimality, so that $P$ is an elementary abelian $p$-group. For any  $c\in Q$ and any $x\in P$ the subgroup $\langle cx, cx^\f,\dots \rangle$ is nilpotent. In particular, this subgroup contains $x^{-1}x^\f=[x,\f]$ and $cx$, so that
 $$
 1=[[x,\f], cx,cx,\dots ]=[[x,\f], c,c,\dots ],
  $$
  whence $[[P,\f], Q]=1$, since the orders of $Q$ and $P$ are coprime. Since $C_G(\f)$ is nilpotent, $Q$ and $C_P(\f)$ also commute. As a result, $Q$ acts trivially on $P=[P,\f]\times C_P(\f)$, a contradiction.
\end{proof}

\section{Ad-nilpotency in the associated Lie algebra} \label{s-ass}

By Theorem~\ref{t-fin} on finite groups the proof of Theorem~\ref{t-2} on profinite groups immediately reduces to the case of pronilpotent groups, and induction on the order $|A|$ of the group of automorphisms will provide a reduction to the case where $A$ has prime order. In this section we consider a certain Lie algebra $L_p(G)$ associated with a finitely generated  pro-$p$ group $G$ admitting an automorphism of prime order satisfying the hypotheses of Theorem~\ref{t-2}. We prove that $L_p(G)$ is generated  by finitely many elements such that all commutators in these elements are ad-nilpotent. The section concludes with a `uniform' version of this result for the associated Lie ring of a finite $p$-group, with essentially the same proof, but with uniform bounds for the ad-nilpotency indices depending on the parameters in the hypotheses.

We begin with some basic definitions.
A series of normal subgroups of a group $G$
\begin{equation}\label{e-ser}
G = G_1\geqslant G_2\geqslant \cdots
\end{equation}
is called a \textit{filtration}
if
\begin{equation}\label{e-fil}
[G_i,G_j] \leqslant G_{i+j}\qquad \text{for all}\quad i, j.
 \end{equation}
 For any filtration \eqref{e-ser} one can define the associated Lie ring $L(G)$ with the additive group
$$
L(G)=\bigoplus _{i}G_i/G_{i+1}.
$$
The Lie product (on the left) is defined for
elements $xG_{i+1}\in L_i$, $yG_{j+1}\in L_j$ via the group commutators (on the right) by
$$
[xG_{i+1},\, yG_{j+1}] = [x, y]G_{i+j+1}\in L_{i+j}
$$
and extended to arbitrary elements of $L(G)$ by linearity. Condition~\eqref{e-fil} ensures that this Lie product is well-defined, and group commutator identities imply that $L(G)$ with these operations is a Lie ring. If all factors $G_i/G_{i+1}$ of a filtration \eqref{e-ser} have prime exponent~$p $,
then $L(G)$ can be viewed as a Lie algebra over the field of $p $ elements $\mathbb{F}_p $.

If all terms of \eqref{e-ser} are invariant under an automorphism $\varphi$ of $G$, then $\varphi$  naturally induces an automorphism of $L(G)$. Note that every homogeneous component $L_{(i)}=G_i/G_{i+1}$ is a $\f$-invariant additive subgroup of $L(G)$. If $\f$ is a coprime automorphism  of a pronilpotent group $G$, then $\f$ induces a coprime automorphism on every factor $G_i/G_{i+1}$ and in particular $L_{(i)}=C_{L_{(i)}}(\f)\oplus [L_{(i)},\f]$ for any $i$ by Lemma~\ref{l-nakr}(c).

In this paper we shall actually use two kinds of associated Lie rings. One of them corresponds to the case where the filtration \eqref{e-ser} is the lower central series, the terms of which are denoted by $\gamma _1(G)=G$ and $\gamma _{i+1}(G)=[\gamma _i(G),G]$.  Note that the corresponding associated Lie ring $L(G)$ is generated by the homogeneous component $L_1=G/\gamma _2(G)$.

 Another example, for a fixed prime number $p $, is the \textit{Zassenhaus $p $-filtration} (also called the \textit{$p $-dimension series}), which is defined by
$$
G_i=\langle g^{p ^k}\mid g\in \gamma _j(G),\;\, jp ^k\geqslant i\rangle .
$$
 The factors of this filtration are elementary abelian $p $-groups, so the corresponding associated Lie ring $D_p (G)$ is a Lie algebra over $\mathbb{F} _p $. We denote  by $L_p (G)$ the subalgebra generated by the first factor $G/D_2$. It is well known that the homogeneous components of $D_p (G)$ of degree $s$ coincide with the homogeneous components of  $ L_p (G)$ for all $s$ that are not divisible by $p $. In particular, $L_p (G)$ is nilpotent if and only if $D_p (G)$ is nilpotent. (Sometimes, the notation $L_p (G)$ is used for $D_p (G)$.)

Recall that an element $a$ of a Lie algebra $L$ is said to be \emph{ad-nilpotent} if there is a number $s$ such that $[l,{}_sa]=0$ for any $l\in L$.

\begin{proposition}
  \label{p-adn}
Suppose that $G$ is a finitely generated pro-$p$ group admitting a coprime automorphism $\f$ of prime order $q$ such that  every subgroup $\langle x^{\langle\f\rangle}\rangle$ generated by the $\langle\f\rangle$-orbit of an element $x\in G$ is nilpotent  and the fixed-point subgroup $C_G(\f)$ is locally nilpotent. Then the Lie algebra $L_p(G)$ is generated by finitely many elements such that all commutators in these elements are ad-nilpotent.
\end{proposition}

\begin{proof}
If elements  $g_1,\dots ,g_d$ (topologically) generate $G$, then their images $\bar g_1,\dots ,\bar g_d$ in $G/D_2(G)$ generate $L_p(G)$. Since $L_p(G)$ is a finitely generated Lie algebra in characteristic $p$, by Lemma~\ref{l-charp} an element $h\in L_p(G)$ is ad-nilpotent if there are numbers $k_1,\dots ,k_n$ such that $[\bar g_i,{}_{k_i}h]=0$ for all $i=1,\dots ,n$.  Therefore in the proof we can safely say that  $h$ is ad-nilpotent if for any $l\in L_p(G)$ there is a number $s$, possibly depending on $l$ (and $h$) such that $[l,{}_sh]=0$. (Although this is actually needed only for $l$ being one of the generators $\bar g_i$, it makes no difference in our arguments.)

\begin{lemma}\label{l-cleft2}
For any $c\in C_G(\f)$ and $g\in G$, the subgroup $\langle c,[g,\f]^{\langle\f\rangle}\rangle$ generated by $c$ and the $\langle\f\rangle$-orbit of $[g,\f]$ is nilpotent.
\end{lemma}

\begin{proof}
  The subgroup $\langle (c[g,\f])^{\langle\f\rangle}\rangle$ generated by the $\langle\f\rangle$-orbit of $c[g,\f]$ is nilpotent by hypothesis. By Lemma~\ref{l-cleft1} this subgroup coincides with the subgroup $\langle c,[g,\f]^{\langle\f\rangle}\rangle$ generated by $c$ and the $\langle\f\rangle$-orbit of $[g,\f]$. Hence the result.
\end{proof}

Although our aim is proving ad-nilpotency of certain elements in $L_p(G)$, we shall perform some calculations over the extended ground field in the Lie algebra $\tilde L=L_p(G)\otimes_{\Z}\Z[\w]$, where $\w$ is a primitive $q$-th root of unity.  The automorphism of $L$ and of $\tilde L$ induced by $\varphi$ is denoted by the same letter. Since $q=|\f|$ is coprime to the characteristic $p$ of the ground field, we have the decomposition $\tilde L=\bigoplus L_i$, where $L_i=\{v\in \tilde L\mid v^\f=\w^iv\}$ are  eigenspaces. In particular, for every $v\in L_(G)$,
\begin{equation}\label{e-eigen}
v=\sum_{i=0}^{q-1}v_i,\qquad \text{where}\quad v_i=\frac{1}{q} \sum_{s=0}^{q-1}\w ^{-is}v^{\f^s};
\end{equation}
here $1/q$ makes sense in characteristic $p$.
This decomposition can also be viewed as a $(\mathbb{Z} /q\mathbb{Z})$-grading of $\tilde L$, since
$$
[L_i,L_j]\subseteq L_{i+j\,({\rm mod}\,q)}.
$$
Note that $L_0=C_{L_p(G)}(\varphi )\otimes _{\mathbb{Z}}\mathbb{Z} [\omega ]$.

We call the intersection $L_{(k)}=L_p(G)\cap D_k/D_{k+1}$ the homogeneous component of weight $k$ of the Lie algebra $L_p(G)$.
Then $[L_{(i)},L_{(j)}]\leq L_{(i+j)}$ for all $i,j$.
We also let $\tilde L _{(k)}=L_{(k)}\otimes _{\mathbb{Z}}\mathbb{Z} [\omega ]$. For clarity we say that elements of  $\tilde L _{(k)}$ or $L _{(k)}$ are \textit{homogeneous of weight $k$}. The component $L_{(1)}$ generates the Lie algebra $L$, and $\tilde L_{(1)}$ generates $\tilde L$. If elements $g_1,\dots ,g_d$ generate the group $G$, then their images $\bar g_1,\dots ,\bar g_d$ in $L_{(1)}=G/D_2$ generate the Lie algebra $L_p(G)$, as well as $\tilde L$ (over the extended ground field).

We set for brevity $\tilde L_{(m)k}= \tilde L_{(m)}\cap L_k$ for any weight $m$. A commutator in homogenous eigenvectors $u_{i}\in L_{(m_i)k_i}$  belongs to $\tilde L_{(m)k}$, where $m$ is the sum of the weights $m_i $ with multiplicities, and $k$ is the modulo $q$ sum of the $k_i$ with multiplicities.
Since $$\tilde L_{(m)}=L_{(m)}+\w L_{(m)}+\w^2L_{(m)}+\dots +\w^{q-2}L_{(m)},$$  any element of $\tilde L_{(m)i}$ is a linear combination of elements of the form $v_i$ in \eqref{e-eigen} for $v\in  L_{(m)}$.
In particular, any homogeneous element $v_0\in \tilde L_{(m)0}=C_{ \tilde L_{(m)}}(\f)$ is a linear combination of the images of elements of $C_{D_m}(\f)$ in $D_m/D_{m+1}$.

We are actually going to prove that  any homogeneous element of $L_p(G)$
is ad-nilpotent, and we approach this goal step-by-step.

\begin{lemma}   \label{l-c0}
Every homogeneous element $\bar c\in C_{L_p(G)}(\f)$ is ad-nilpotent in $\tilde L$.
\end{lemma}

\begin{proof}
Since  $\tilde L$ is a direct sum of homogeneous components, it is sufficient to prove that for any $\bar w$ in a homogeneous component $\tilde L_{(a)}$ we have $[\bar w,{}_s\bar c]=0$ for some $s$. Since $\tilde L_{(a)}=L_{(a)}\otimes _\Z \Z[\w]$, any such $\bar w$ is a finite linear combination of homogeneous elements in  $L_{(a)}$. Thus, it suffices to show that   for any $\bar u\in  L_{(a)}$ we have $[\bar u,{}_s\bar c]=0$ for some $s$.

If $\bar c\in D_b/D_{b+1}$ is homogeneous of weight $b$, let $c\in C_G(\f)\cap D_b$ be a pre-image of $\bar c$. For any $\bar z\in C_{L_{(a)}}(\f)$ we obviously have $[\bar z,{}_s\bar c]=0$ for some $s$, since the subgroup $\langle z,c\rangle$ generated by the pre-images of $\bar z$ and $\bar c$ in $C_G(\f)$ is nilpotent by hypothesis.
Any element of $[L_{(a)},\f ]$ has the form $[\bar g,\f]$ for some $\bar g\in L_{(a)}$. Let $g\in D_k$ be a pre-image of $\bar g$ in~$D_a$; then  $[\bar g,\f]=\overline{[g,\f]}$ is the image of $[g,\f]$.
By Lemma~\ref{l-cleft2} we have $[[g,\f], {}_s c]=1$ for some $s$. Hence,
 $[[ \bar g,\f], {}_s\bar c]=0$ in $L_p(G)$.
Since any element of $L_{(a)}$ has the form $[\bar g,\f]+\bar z$ for $\bar z\in C_{L_{(a)}}(\f)$ and $\bar g\in L_{(a)}$,  we obtain that $[[\bar g,\f]+\bar z,{}_s\bar c]=0$ for some $s$, as required.
\end{proof}

\begin{lemma}\label{l-adn1}
Every homogeneous element $c\in C_{\tilde L}(\f)$ is ad-nilpotent.
\end{lemma}

\begin{proof}
We have $c= \bar c_0+ \w c_1+\w^2 \bar c_2\cdots +w^{q-2} \bar c_{q-2}$ for  some homogeneous elements $\bar c_i\in C_{L_p(G)}(\f)$, which are the images of some elements $c_i\in C_{D_a}(\f)$ in $D_a/D_{a+1}$, where $a$ is the weight of $c$.  Every commutator in the elements $\bar c_0, \w  \bar c_1,\dots ,w^{q-2} \bar c_{q-2}$ is ad-nilpotent
by Lemma~\ref{l-c0}, since it has the form $\w ^t \bar \ka $ for some $t$ and a homogeneous element $\bar \ka \in C_{L_p(G)}(\f )$ (the corresponding commutator in the $\bar c_i$). The subgroup $\langle c_0,\dots , c_{q-2}\rangle$ is nilpotent;
hence the Lie subalgebra $\langle\bar c_0, \w  \bar c_1,\dots ,w^{q-2} \bar c_{q-2}\rangle$ is also nilpotent. Hence there are only finitely many commutators in the elements $\bar c_0, \w  \bar c_1,\dots ,w^{q-2} \bar c_{q-2}$, and so there is $n$ such that all these commutators are ad-nilpotent of index $n$. By Lemma~\ref{l-be} every element of this Lie subalgebra is ad-nilpotent, including $c=\bar c_0+ \w c_1+\w^2 \bar c_2\cdots +w^{q-2} \bar c_{q-2}$.
\end{proof}

\begin{lemma}\label{l-adn2}
Every homogeneous element $[\bar g,\f]\in [L_p(G), \f]$ is ad-nilpotent in $\tilde L$.
\end{lemma}

\begin{proof}
Let $a$ be the weight of $[\bar g,\f]$, so that $[\bar g,\f]\in [L_{(a)}, \f]$.
We have the eigenvector decomposition
\begin{equation}\label{e-decom}
  [\bar g,\f]=b_1+\dots +b_{q-1},\qquad\text{where}\quad b_i=\frac{1}{q} \sum_{s=0}^{q-1}\w ^{-is}[\bar g,\f]^{\f^s}\in \tilde L_i.
\end{equation}
(Here, $b_0=(1/q) \sum_{s=0}^{q-1}[\bar g,\f]^{\f^s}=0$.)
 We claim that any commutator in the $b_i$ is ad-nilpotent. Any such commutator $\ka _j$ is also an eigenvector corresponding to $\w ^j$, where $j$ is the mod $q$ sum of the indices of all of its entries $b_i$ (with multiplicities). Since $\tilde L=\tilde L_0+\dots +\tilde L_{q-1}$, it is sufficient to prove that for any homogeneous eigenvector $\lambda _u\in \tilde L_u$ we have $[\lambda _u,{}_s\ka_j]=0$ for sufficiently large $s$. By Lemma~\ref{l-adn1} we can assume that $j\ne 0$.
 Since $q$ is a prime, there is $v\in\{0,1,\dots ,q-1\}$ such that  $u+vj\equiv 0\,({\rm mod}\, q)$.
  Replacing $\lambda _u$ with the
initial segment of length $1+vj$, we can assume that $\lambda _u=\lambda _0\in \tilde L_0$.  Such $\lambda _0$ is a linear combination of homogeneous elements of $C_{L_p(G)}(\f)$ of some weight $m$. Therefore it is sufficient to show that  for $\bar c\in C_{L_{(m)}}(\f)$ we have $[\bar c,{}_s\ka_j]=0$ for some~$s$.

 Let $c\in C_{D_m}(\f)$ be a pre-image of $\bar c$, and $g\in D_a$ a pre-image of $\bar g$.
By Lemma~\ref{l-cleft2} the subgroup $\langle c, [g,\f]^{\langle\f\rangle}\rangle$ is nilpotent of some class $n$. Therefore,
$$
  \big[c,\, [g,\f]^{\f^{s_1}}, [g,\f]^{\f^{s_2}},\dots , [g,\f]^{\f^{s_n}}\big]=1
$$
for every $n$-tuple $(s_1,\dots s_n)$.  Hence the Lie algebra commutator
\begin{equation}\label{e-3}
\big[\bar c, \,[\bar g,\f]^{\f^{s_1}}, [\bar g,\f]^{\f^{s_2}},\dots , [\bar g,\f]^{\f^{s_n}}\big]
\end{equation}
 is equal to~0.

 We can now prove that $[\bar c,{}_n\ka_j]=0$. We replace all the entries $b_i$ in $\ka _j$ with their expressions \eqref{e-decom} as linear combinations of the $[\bar g,\f]^{\f^{i}}$, $i=0,1,\dots ,q-1$. Repeatedly using the Jacobi identity $[a,[b,c]]=[[a,b],c]-[[a,c],b]$  we obtain a linear combination of commutators with initial segments of the form \eqref{e-3}, which are all trivial.

 We finally prove that $[\bar g,\f]$ is ad-nilpotent.  We showed above that any commutator in the $b_i$ is ad-nilpotent. The subalgebra $\langle b_1,\dots ,b_{q-1}\rangle$ is nilpotent, since it is contained in (actually, is equal to) the subalgebra $\langle [\bar g,\f],\dots ,[\bar g,\f]^{\f^{q-1}}\rangle$, which is nilpotent, since the subgroup $\langle [ g,\f],\dots ,[ g,\f]^{\f^{q-1}}\rangle$ generated by pre-images in $D_a$ is nilpotent by hypothesis.   Hence there are only finitely many commutators in the $b_i$, and so there is $n$ such that all these commutators are ad-nilpotent of index $n$. Then by Lemma~\ref{l-be} every element of the subalgebra $\langle b_1,\dots ,b_{q-1}\rangle$ is ad-nilpotent, including
  $[\bar g,\f]=b_1+\dots +b_{q-1}$.
\end{proof}

We can now finish the proof of Proposition~\ref{p-adn}.

Recall that we want to prove that the Lie algebra $L_p(G)$ is generated by finitely many elements such that all commutators in these elements are ad-nilpotent.
The Lie algebra $L_p(G)$ is clearly  generated by finitely many homogeneous elements of weight 1, since $G$ is finitely generated.  We claim that actually every homogeneous element $\bar g\in L_{(a)}$  for any weight $a$ is ad-nilpotent.
Let $M=\langle \bar g, \bar g^\f,\dots,\bar g^{\f^{q-1}}\rangle$ be the minimal $\f$-invariant Lie subalgebra containing $\bar g$. For a pre-image $g\in D_a$ of $\bar g$, the subgroup $\langle  g,  g^\f,\dots, g^{\f^{q-1}}\rangle$ is nilpotent of some class $n$; hence $M$ is also nilpotent of class $n$. For $s=1,\dots ,n$ let $M_{(s)}$ be the additive subgroup generated by all commutators of weight $s$ in the  $\bar g^{\f^{i}}$. This is a $\f$-invariant additive subgroup of the homogeneous component $L_{(as)}$ of weight $as$. Note that $M_{(s)}=C_{M_{(s)}}(\f)\oplus  [M_{(s)},\f]$.

We now prove by reverse induction on $i$ that every element $\bar u$ of $M_{(i)}$ is ad-nilpotent. At the base of induction, for $i=n$, we have $\bar u=\bar z+\bar w$ for $\bar z\in C_{M_{(n)}}(\f)$ and $\bar w\in  [M_{(n)},\f]$. We know that $\bar z$ is ad-nilpotent by Lemma~\ref{l-c0}, and $\bar w$ is ad-nilpotent by Lemma~\ref{l-adn2}. Since they commute, the sum is also ad-nilpotent: for any $l\in L_p(G)$ we have $[l,{}_s\bar z]=0$ for some~$s$, and $[l,{}_t\bar w]=0$ for some $t$, and then $[l,{}_{s+t}(\bar z+\bar w)]=0$.

For  $\bar u\in M_{(i)}$ we have $\bar u=\bar z+\bar w$ for $\bar z\in C_{M_{(i)}}(\f)$ and $\bar w\in  [M_{(i)},\f]$. We know that $\bar z$ is ad-nilpotent by Lemma~\ref{l-c0}, and $\bar w$ is ad-nilpotent by Lemma~\ref{l-adn2}. Every commutator of weight $\geq 2$ in $\bar z$ and $\bar w$ belongs to $M_{(j)}$ for some $j>i$ and therefore is also ad-nilpotent by the induction hypothesis. The Lie subalgebra $\langle  \bar z ,\bar w\rangle$ is nilpotent. Hence there are only finitely many commutators in the elements $ \bar z ,\bar w$, and so there is $n$ such that all these commutators are ad-nilpotent of index $n$. Hence by Lemma~\ref{l-be} every element of this subalgebra is ad-nilpotent, including $\bar u$. At the end of this induction, when  $i=1$, we obtain that $\bar g$ is ad-nilpotent.
\end{proof}

We now state a `uniform' version of Proposition~\ref{p-adn}, in which $L(G)$ is the associated Lie ring  constructed from the lower central series.

\begin{proposition}
  \label{p-adn-u}
  Let $d,k,l$ be positive integers.
Suppose that $G$ is a finite $p$-group generated by $d$ elements. Suppose that $G$ admits an automorphism $\f$ of prime order $q\ne p$ such that  every subgroup $\langle x^{\langle\f\rangle}\rangle$ generated by the $\langle\f\rangle$-orbit of an element $x\in G$ is nilpotent of class $l$ and the fixed-point subgroup $C_G(\f)$ is nilpotent of class $k$.
Then the associated Lie ring $L(G)$ is generated by $d$ elements such that all commutators $\ka$ in these elements are ad-nilpotent of  uniformly bounded index, that is, $[L(G),{}_{f}\ka ]=0$ for some $f=f(q,k,l)$ depending only on $q,k,l$.
\end{proposition}

The proof follows closely the proof of Proposition~\ref{p-adn}, and we only indicate the adjustments which ensure uniform bounds for ad-nilpotency indices. We use the same notation for the Lie ring $\tilde L=L(G)\otimes_\Z \Z[\w]$ and elements of homogeneous components as in the corresponding lemmas in the proof of Proposition~\ref{p-adn}, with $L(G)$ replacing $L_p(G)$ throughout. In particular, we have a similar decomposition into analogues of eigenspaces $\tilde L=\bigoplus L_i$, where $L_i=\{v\in \tilde L\mid v^\f=\w^iv\}$, because the additive group of $L(G)$ has unique $q$-th roots of all elements, since it is an abelian $p$-group for $p\ne q$.

\begin{lemma}\label{l-cleft2-u}
For any $c\in C_G(\f)$ and $g\in G$, the subgroup $\langle c,[g,\f]^{\langle\f\rangle}\rangle$ generated by $c$ and the $\langle\f\rangle$-orbit of $[g,\f]$ is nilpotent of class $l$.
\end{lemma}

\begin{proof}
The subgroup $\langle (c[g,\f])^{\langle\f\rangle}\rangle$ generated by the $\langle\f\rangle$-orbit of $c[g,\f]$ is nilpotent of class~$l$. By Lemma~\ref{l-cleft1} this subgroup coincides with the subgroup $\langle c,[g,\f]^{\langle\f\rangle}\rangle$ generated by $c$ and the $\langle\f\rangle$-orbit of $[g,\f]$. Hence the result.
\end{proof}

\begin{lemma}   \label{l-c0-u}
Every homogeneous element $\bar c\in C_{L(G)}(\f)$ is ad-nilpotent in $\tilde L$ of $(k,l)$-bounded index, that is, $[\tilde L,{}_{f_1}\bar c]=0$ for some $f_1=f_1(k,l)$ depending only on $k$ and $l$.
 \end{lemma}

\begin{proof}
We repeat the proof of the corresponding Lemma~\ref{l-c0} with the following modifications. For any $\bar z\in C_{L_{(a)}}(\f)$ we have $[\bar z,{}_{k}\bar c]=0$, since the subgroup $\langle z,c\rangle$ is nilpotent of class~$k$. Further, by Lemma~\ref{l-cleft2-u} we have $[[g,\f], {}_l c]=1$, whence
 $[[ \bar g,\f], {}_l\bar c]=0$ in $L(G)$. Then for $s=\max\{k,l\}$ we obtain $[[\bar g,\f]+\bar z,{}_s\bar c]=0$.
\end{proof}

\begin{lemma}\label{l-adn1-u}
Every homogeneous element $c\in C_{\tilde L}(\f)$ is ad-nilpotent of $(q,k,l)$-bounded index, that is, $[\tilde L,{}_{f_2}c ]=0$ for some $f_2=f_2(q,k,l)$ depending only on $q,k,l$.
\end{lemma}

\begin{proof} We repeat the proof of the corresponding Lemma~\ref{l-adn1} with the following modifications.
 Every commutator in the elements $\bar c_0, \w  \bar c_1,\dots ,w^{q-2} \bar c_{q-2}$ (including these elements) is ad-nilpotent of index $f_1(k,l)$
by Lemma~\ref{l-c0-u}. The subgroup $\langle c_0,\dots , c_{q-2}\rangle$ is nilpotent of class $k$; hence  the Lie subring $\langle\bar c_0, \w  \bar c_1,\dots ,w^{q-2} \bar c_{q-2}\rangle$ is also nilpotent of class $k$. By Lemma~\ref{l-be} every element of this Lie subring is ad-nilpotent of  $(q,k,l)$-bounded index  $f_2(q,k,l)$, including $c=\bar c_0+ \w c_1+\w^2 \bar c_2\cdots +w^{q-2} \bar c_{q-2}$.
\end{proof}

\begin{lemma}\label{l-adn2-u}
Every homogeneous element $[\bar g,\f]\in [L(G), \f]$ is ad-nilpotent in $\tilde L$ of $(q,k,l)$-bounded index, that is, $[\tilde L,{}_{f_3}[\bar g,\f]]=0$ for some $f_3=f_3(q,k,l)$ depending only on $q,k,l$.
\end{lemma}

\begin{proof}
We repeat the proof of the corresponding Lemma~\ref{l-adn2} with the following modifications.
 By Lemma~\ref{l-adn1-u} any $\ka _0$ is ad-nilpotent of  index  $f_2(q,k,l)$; hence we can assume that $j\ne 0$. Replacing $\lambda _u$ with the
initial segment of length $1+vj$, where  $u+vj\equiv 0\,({\rm mod}\, q)$ for some $v\leq q-1$, we can assume that $\lambda _u=\lambda _0\in \tilde L_0$.
By Lemma~\ref{l-cleft2-u} the subgroup $\langle c, [g,\f]^{\langle\f\rangle}\rangle$ is nilpotent of class $l$. Hence,
$$
  \big[c,\, [g,\f]^{\f^{s_1}}, [g,\f]^{\f^{s_2}},\dots , [g,\f]^{\f^{s_l}}\big]=1
$$
for every $l$-tuple $(s_1,\dots s_l)$.  Hence the Lie ring commutator
$$
\big[\bar c, \,[\bar g,\f]^{\f^{s_1}}, [\bar g,\f]^{\f^{s_2}},\dots , [\bar g,\f]^{\f^{s_l}}\big]
$$
 is equal to~0. It follows that $[\bar c,{}_l\ka_j]=0$. Then also $[\bar \lambda_u,{}_{l+v}\ka_j]=0$. Thus,  $\ka_j$ is ad-nilpotent of index $l+q-1$ if $j\ne 0$, and of index $f_2(q,k,l)$ for $j=0$.

The subring $\langle b_1,\dots ,b_{q-1}\rangle$ is nilpotent of class $l$, since it is contained in  the subring $\langle [\bar g,\f],\dots ,[\bar g,\f]^{\f^{q-1}}\rangle$, which is nilpotent of class $l$, since the subgroup $\langle [ g,\f],\dots ,[ g,\f]^{\f^{q-1}}\rangle$ is nilpotent of class $l$.
Since any commutator in the $b_i$ is ad-nilpotent of index $\max\{l+q-1,f_2(q,k,l)\}$, by Lemma~\ref{l-be} every element of the subring $\langle b_1,\dots ,b_{q-1}\rangle$ is ad-nilpotent of $(q,k,l)$-bounded index $f_3(q,k,l)$, including
  $[\bar g,\f]=b_1+\dots +b_{q-1}$.
\end{proof}

We now finish the proof of Proposition~\ref{p-adn-u}.

\begin{proof}[Proof of Proposition~\ref{p-adn-u}]
We repeat the proof of Proposition~\ref{p-adn} with the following modifications.
 We claim that actually every homogeneous element $\bar g\in L_{(a)}$  for any weight $a$ is ad-nilpotent of $(q,k,l)$-bounded index.
The Lie subring $M=\langle \bar g, \bar g^\f,\dots,\bar g^{\f^{q-1}}\rangle$
is nilpotent of class $l$.
We prove by reverse induction on $i$ that every element $\bar u$ of $M_{(i)}$ is ad-nilpotent of  $(q,k,l)$-bounded index. For $i=n$, we have $\bar u=\bar z+\bar w$ for $\bar z\in C_{M_{(n)}}(\f)$ and $\bar w\in  [M_{(n)},\f]$, where $\bar z$ is ad-nilpotent of  $(k,l)$-bounded index by Lemma~\ref{l-c0-u}, and $\bar w$ is ad-nilpotent of $(q,k,l)$-bounded index by Lemma~\ref{l-adn2-u}. Since they commute, the sum is also ad-nilpotent of  $(q,k,l)$-bounded index.

For  $\bar u\in M_{(i)}$ we have $\bar u=\bar z+\bar w$ for $\bar z\in C_{M_{(i)}}(\f)$ and $\bar w\in  [M_{(i)},\f]$. We know that $\bar z$ is ad-nilpotent of  $(k,l)$-bounded index by Lemma~\ref{l-c0-u}, and $\bar w$ is ad-nilpotent of  $(q,k,l)$-bounded index by Lemma~\ref{l-adn2-u}. Every commutator of weight $\geq 2$ in $\bar z$ and $\bar w$ belongs to $M_{(j)}$ for some $j>i$ and therefore is also ad-nilpotent of  $(q,k,l)$-bounded index by the induction hypothesis. The Lie subring $\langle  \bar z ,\bar w\rangle$ is nilpotent of class $l$. Hence by Lemma~\ref{l-be} every element of this subring is ad-nilpotent  of $(q,k,l)$-bounded index, including $\bar u$. At the end of this induction, when  $i=1$, we obtain that $\bar g$ is ad-nilpotent of $(q,k,l)$-bounded index.
\end{proof}

\section{Pro-$p$ groups}

In this section we prove Theorem~\ref{t-2} for pro-$p$ groups.

\begin{theorem}\label{t-pro-p}
Suppose that a pro-$p$ group $G$ admits a soluble finite group $A$ of coprime automorphisms such that for every $\f\in A$  every subgroup $\langle x^{\langle\f\rangle}\rangle$ generated by the $\langle\f\rangle$-orbit of an element $x\in G$ is nilpotent. If the fixed-point subgroup $C_G(A)$ is locally nilpotent, then $G$ is  locally nilpotent.
\end{theorem}

\begin{proof}
  First we note that it is sufficient to prove the theorem for the case where $A$ has prime order.
  Indeed,  since $A$ is soluble, there is a normal subgroup $A_1$ of prime index $q$. Let $\f\in A\setminus A_1$. Then  $C_G(A_1)$ admits the cyclic group of automorphisms of order $q$ induced by $\langle\f\rangle$, which satisfies on $C_G(A_1)$ the same hypotheses as $A$. Assuming that the theorem is proved in the case where $A$ has prime order, we obtain that $C_G(A_1)$ is nilpotent. Now the group $G$ admits the group of automorphisms $A_1$ of smaller order satisfying the same hypotheses, and therefore the proof concludes by the obvious induction on  $|A|$.   Thus, from now on we assume that $G=\langle\f\rangle$ is cyclic of prime order $q$.

Any finite subset $S\subseteq G$ is contained in the $\varphi$-invariant finite set $S^{\langle\varphi \rangle}=\{s^{\varphi ^k}\mid s\in S,\; k=0,1,\dots ,|\f |-1\}$. Therefore we can assume that $G$ is (topologically) finitely generated, and then we need to prove that $G$ is nilpotent.

Thus, the group $G$ and its automorphism $\f$ now satisfy the hypotheses of Proposition~\ref{p-adn}.

 Let $L_p (G)$  be the associated Lie algebra of  $G$ over $\mathbb{F} _p $ constructed with respect to the Zassenhaus $p $-filtration of $G$ as described in \S\,\ref{s-ass}.
 First we shall prove that $L_p (G)$ is nilpotent using Zelmanov's Theorem~\ref{t-z}. This will imply that $G$ is a $p$-adic analytic group by the theory of Lazard \cite{laz} and Lubotzky--Mann \cite{lu-ma}. In particular, then $G$ is a linear group and has finite rank, whence $C_G(\f)$ is nilpotent. We then proceed by induction on the nilpotency class of $C_G(\f)$ to show that $G$ is soluble. At the basis of this induction we have the case $C_G(\f)=1$, which implies  that $G$ is nilpotent of class $h(q)$ by Higman's theorem~\cite{hig}. Induction on the derived length then finishes the proof.

\begin{proposition}\label{p-lpg}
The Lie algebra $L_p (G)$ is nilpotent.
\end{proposition}

\begin{proof}
By Zelmanov's Theorem~\ref{t-z}, $L_p(G)$ is nilpotent if it satisfies a polynomial identity and is generated by finitely many elements such that all commutators in these elements are ad-nilpotent. The second condition is satisfied by Proposition~\ref{p-adn}. It remains to show that $L_p(G)$ satisfies a polynomial identity.  Recall that a group $G$ is said to satisfy a \textit{coset identity} if there is a group word $w(x_1,\dots ,x_m)$ and cosets $a_1H,\dots ,a_mH$ of a subgroup $H\leqslant G$ of finite index such that $w(a_1h_1,\dots ,a_mh_m)=1$ for any $h_i\in H$. Wilson and Zelmanov \cite{wi-ze} proved that if a group $G$ satisfies a coset identity, then the Lie algebra $L_p (G)$ constructed with respect to the Zassenhaus $p $-filtration satisfies a polynomial identity.

\begin{lemma}\label{l-pi}
The Lie algebra $L_p (G)$ satisfies a polynomial identity.
\end{lemma}

\begin{proof}
Since $ C_G(\varphi )$ is locally nilpotent, it follows that $ C_G(\varphi )$ satisfies a coset identity on cosets of an open subgroup of $C_G(\varphi )$. For example, in the group $ C_G(\varphi )\times C_G(\varphi )$ the subsets
$$
E_i=\{ (x,y)\in C_G(\varphi )\times C_G(\varphi ) \mid [x,{}_iy]=1\}
$$
are closed in the product topology, and
$$
C_G(\varphi )\times C_G(\varphi )=\bigcup_{i=1}^{\infty} E_i.
$$
Hence by the Baire category theorem \cite[Theorem~34]{kel}, one of these subsets $E_n$ contains an open subset of  $ C_G(\varphi )\times C_G(\varphi )$, which means that there are cosets $x_0K_1$, $y_0K_2$ of open subgroups $K_1,K_2\leqslant C_G(\varphi )$ such that $[x,{}_ny]=1$ for all $x\in x_0K_1$ and $y\in y_0K_2$, and therefore for all $x\in x_0(K_1\cap K_2)$ and $y\in y_0(K_1\cap K_2)$. Thus, $ C_G(\varphi )$ satisfies a coset identity.

The intersections $C_i=C_G(\varphi )\cap D_i$ with the terms $D_i$ of the Zassenhaus $p $-filtration for $G$ form a filtration of $C_G(\varphi )$, since $[C_i, C_j]\leqslant C_{i+j}$. Let $L_c$ be the Lie algebra  constructed for $C_G(\varphi )$ with respect to the filtration $\{ C_i \}$ in the same way as $D(G)$. Since $\varphi$ is a coprime automorphism, the fixed-point subalgebra $ C_{D(G)}(\varphi )$ is isomorphic to $L_c$. We apply a version of the Wilson--Zelmanov result \cite[Theorem~1]{wi-ze}, by which a coset identity on $C_G(\varphi )$ implies that $L_c$ satisfies some polynomial identity. Indeed, the proof of Theorem~1 in \cite{wi-ze} only uses the filtration properties $[F_i,F_j]\leqslant F_{i+j}$ and $F_i^p\leq F_{ip}$ for showing that the homogeneous Lie polynomial constructed from a coset identity on a group $F$ is an  identity of the  Lie algebra constructed with respect to the filtration $\{F_i\}$, and these properties are also satisfied by the filtration $\{C_i\}$.

Thus, the fixed-point subalgebra $ C_{D(G)}(\varphi )$ satisfies a polynomial identity. Hence the Lie algebra $D(G)$ also satisfies a polynomial identity by the Bahturin--Zaicev Theorem~\ref{t-bz}, and therefore so does $L_p(G)$.
\end{proof}

We are now ready to prove that $L_p(G)$ is nilpotent.
By Lemma~\ref{l-pi} the Lie algebra $L_p(G)$ satisfies a polynomial identity, and by Proposition~\ref{p-adn} it is generated by finitely many elements such that all commutators in these elements are ad-nilpotent. Hence $L_p(G)$ is nilpotent by Zelmanov's Theorem~\ref{t-z}.
\end{proof}

We now return to the proof of Theorem~\ref{t-pro-p}, that is, of the nilpotency of~$G$. The nilpotency of the Lie algebra $L_p (G)$  of the finitely generated pro-$p $ group $G$  established in Proposition~\ref{p-lpg}  means that $G$ is a $p $-adic analytic group and therefore a linear group. This result goes back to Lazard~\cite{laz}; see also \cite[Corollary~D]{sha}. By the Lubotzky--Mann theorem \cite{lu-ma}, being a $p $-adic analytic group is equivalent to being a pro-$p $ group of finite rank. In particular, all closed subgroups are finitely (topologically) generated, and therefore $C_G(\f)$ is nilpotent of some class $k$.   We now show that $G$ is soluble.

\begin{proposition}
  \label{p-pro-p-sol}
  The group $G$ is soluble.
\end{proposition}

\begin{proof}
  We proceed by induction on the nilpotency class $k$ of $C_G(\f)$. At the basis is the case of $C_G(\f)=1$, which formally corresponds to $k=0$. Then $G$ is nilpotent of class $h(q)$ by Higman's theorem \cite{hig} applied to the finite quotients of $G$ by $\f$-invariant normal subgroups. Now let $k\geq 1$.

  \begin{lemma}
    \label{l-left}
    Every element of the centre $Z(C_G(\f))$ is a left Engel element of $G$.
  \end{lemma}

  \begin{proof}
    Let $z\in Z(C_G(\f))$ and $g\in G$. By Lemma~\ref{l-coprp} there are elements $b\in G$ and $c\in C_G(\f))$ such that $g=c[b,\f]$. Then $[g,z]=[[b,\f],z]$. The group $\langle z, [b,\f]^{\langle\f\rangle}\rangle$ is nilpotent of some class $s$ by Lemma~\ref{l-cleft2}.  Then $[g,{}_sz]=[[b,\f],{}_sz]=1$.
  \end{proof}

  By Gruenberg's theorem \cite{gru} left Engel elements of a linear group belong to its abstract Hirsch--Plotkin radical $R(G)$, which is  locally nilpotent. Since the rank of $G$ is finite, $R(G)$ is a closed subgroup. Indeed, the Frattini quotient of the closure $\overline{R(G)}/\Phi(\overline{R(G)})$ is finite. By the Burnside basis theorem any set of topological generators contains a finite subset that topologically generates $K$. Hence we can choose finitely many elements in $R(G)$ that topologically generate $\overline{R(G)}$. The abstract subgroup generated by these elements is nilpotent and so $\overline{R(G)}$ is nilpotent and $\overline{R(G)}=R(G)$.

  Since $Z(C_G(\f))\leq R(G)$ by Lemma~\ref{l-left}, the nilpotency class of $C_{G/R(G)}(\f)$ is smaller than~$k$. By the induction hypothesis, $G/R(G)$ is soluble, and therefore $G$ is soluble.
\end{proof}

We now finish the proof of the nilpotency of $G$ by induction on its derived length. The main case is when $G$ is metabelian.

\begin{lemma}
  \label{l-metab}
  If $G$ is metabelian, then it is nilpotent. \end{lemma}

\begin{proof}
First we show that $[G,\f]$ is nilpotent. We claim that every element $c\in C_G(\f)\cap [G,\f]$ is right Engel in $[G,\f]$. Indeed, by Lemma~\ref{l-coprp} every element of $[G,\f]$ has the form $c_1[u,\f]$ for some $c_1\in C_G(\f)\cap [G,\f]$. By Lemma~\ref{l-cleft2} we have $[c,{}_s[u,\f]]=1$ for some $s$. Since both $c$ and $c_1$ belong to the abelian  derived subgroup $[G,\f ]'$, we have  $[c,{}_sc_1[u,\f]]= [c,{}_s[u,\f]]=1$. Thus, $c$ is right Engel in $[G,\f]$.  By Gruenberg's theorem \cite{gru} right Engel elements of a linear group belong to its hypercentre, so that  $C_G(\f)\cap [G,\f]\leq Z_\infty ([G,\f])$. Hence $\f$ has no nontrivial fixed point on  $[G,\f]/Z_\infty ([G,\f])$, whence $[G,\f]/Z_\infty ([G,\f])$ is nilpotent of class $h(q)$ by Higman's theorem \cite{hig}. This means that $[G,\f]= Z_\infty ([G,\f])$, and so $[G,\f]$ is locally nilpotent \cite[12.2.4]{rob}, and even nilpotent, since it has finite rank.

We now show that every element $c\in C_G(\f)$ is left Engel. For any $g\in G$, by Lemma~\ref{l-coprp} we have $[g,c]=c_1[u,\f]$, where $c_1\in C_G(\f)\cap G'$. Then $[g,{}_{s+1}c]=[c_1,{}_sc][[u,\f],{}_sc]$. By Lemma~\ref{l-cleft2} we have $[[u,\f],{}_sc]=1$ for large enough $s$, and $[c_1,{}_sc]=1$ for large enough $s$ since $C_G(\f)$ is nilpotent. Thus, $c$ is a left Engel element.  By Gruenberg's theorem \cite{gru} left Engel elements of a linear group belong to its Hirsch--Plotkin radical $R(G)$, which is closed and nilpotent as shown above in the proof of Lemma~\ref{l-left}.

Now $G=[G,\f]R(G)$ is the product of two nilpotent normal subgroups and therefore is nilpotent.
\end{proof}

We can now finish the proof of Theorem~\ref{t-pro-p} by induction on the derived length of $G$ because $G$ is soluble by Proposition~\ref{p-pro-p-sol}. Note that in a finitely generated pro-$p$ group the derived subgroup is closed by a theorem of Serre \cite[\S\,4.2]{ser}. Since $G$ has finite rank, the derived subgroup $G'$ is also finitely generated. By the induction hypothesis, $G'$ is nilpotent. Since $G/G''$ is nilpotent by Lemma~\ref{l-metab}, the group $G$ is nilpotent by Hall's theorem \cite{hal}.
\end{proof}

\section{Uniform bound for the nilpotency class}

In Theorem~\ref{t-2}, we need to prove that if a finitely generated profinite group $G$ admits a soluble finite group $A$ of coprime automorphisms such that for every $\f\in A$ all subgroups $\langle x^{\langle\f\rangle}\rangle$ are nilpotent and $C_G(A)$ is locally nilpotent, then $G$ is nilpotent. Theorem~\ref{t-fin} implies that $G$ is pronilpotent, and every Sylow $p $-subgroup of $G$ is nilpotent by Theorem~\ref{t-pro-p}. But this does not yet imply the nilpotency of $G$,
since the nilpotency class given by Theorem~\ref{t-pro-p} depends on $p$.
What we need is a uniform bound for the nilpotency classes of Sylow $p$-subgroups of $G$. We achieve such a uniform bound in this section by proving that there is a number $N$ such that for $p>N$ the Sylow $p$-subgroups of $G$ are nilpotent of class bounded in terms of certain parameters of $G$ and $A$ only, independently of $p$. Together with the nilpotency of the Sylow $p$-subgroups for all $p\leq N$, this will complete the proof of the nilpotency of~$G$.

One of the parameters in question is the finite number of generators of $G$ denoted  by~$d$. Clearly, every Sylow $p $-subgroup of $G$ can also be generated by $d$ elements. Another parameter will be the nilpotency class of $C_P(A)$ for the Sylow $p$-subgroups $P$ for all $p>N$.

\begin{lemma}\label{l-uni-e}
There are positive integers $n$ and $N_1$ such that  for every $p >N_1$, in the Sylow $p $-subgroup $P$ of $G$,  the fixed-point subgroup $C_P(A)$ is an $n$-Engel group.
\end{lemma}

\begin{proof}
 In the group $ C_G(A )\times C_G(A )$,
the subsets
$$
E_i=\{ (x,y)\in C_G(A )\times C_G(A ) \mid [x,{}_iy]=1\}
$$
are closed in the product topology, and
$$
C_G(A )\times C_G(A )=\bigcup_{i=1}^{\infty} E_i
$$
because $C_G(A)$ is locally nilpotent. Hence by the Baire category theorem \cite[Theorem~34]{kel}, one of these subsets $E_n$ contains an open subset of  $ C_G(A )\times C_G(A )$. This means that there are cosets $x_0K_1$, $y_0K_2$ of open subgroups $K_1,K_2\leqslant C_G(A )$ such that $[x,{}_ny]=1$ for all $x\in x_0K_1$ and $y\in y_0K_2$, and therefore for all $x\in x_0K$ and $y\in y_0K$, where $K=K_1\cap K_2$. Clearly,  $C_G(A)=\prod_{p}C_{G_p}(A)$ is the Cartesian product of the $C_{G_p}(A)$ over Sylow $p$-subgroups $G_p$ of~$G$. Since $K$ has finite index in $C_G(A)$, there is a positive integer $N_1$ such that $K$ contains $C_{G_p}(A)$ for all $p>N_1$. The coset representatives $x_0,y_0$ then can be chosen in the direct product $\prod_{p\leq N_1}C_{G_p}(A)$. Therefore the projection of the relation $[x,{}_ny]=1$ for  $x\in x_0K$ and $y\in y_0K$ onto $C_{G_p}(A)$ for $p>N_1$ means that $[x,{}_ny]=1$ for all $x,y\in C_{G_p}(A)$.
\end{proof}

\begin{lemma}\label{l-uni-c}
 There are positive integers $k$ and $N_2$ such that for every $p >N_2$, in the Sylow $p $-subgroup $P$ of $G$, the fixed-point subgroup $C_P(A )$ is nilpotent of class $k$.
 \end{lemma}

\begin{proof}
By Lemma~\ref{l-uni-e}, for $p >N_1$ in the Sylow $p $-subgroup $P$ the subgroup  $C_P(A )$ is an $n$-Engel group. By a theorem of Burns and Medvedev \cite{bu-me}, then  $C_P(A )$ has a normal subgroup $N_p $ of exponent $e(n)$ such that the quotient $C_P(A )/N_p $ is nilpotent of class $k(n)$, for some numbers $e(n)$ and $k(n)$ depending only on $n$. Clearly, $N_p =1$ for all large enough primes $p >N_2\geqslant N_1$. Thus, for every prime $p >N_2$ the subgroup $C_P(A )$ is nilpotent of class $k=k(n)$.
\end{proof}

One more parameter that will be used for uniformly bounding the nilpotency class of
the Sylow $p$-subgroups for all large enough primes $p$ is the uniform bound for the nilpotency classes of the subgroups $\langle x^{\langle\f\rangle}\rangle$ given by the following lemma.

\begin{lemma}\label{l-uni-z}
 There are positive integers $l$ and $N_3$ such that for every $p >N_3$, in the Sylow $p$-subgroup $P$ of $G$, for every $\f\in A$ and every $x\in P$ the subgroup  $\langle x^{\langle\f\rangle}\rangle$ generated by the $\langle\f\rangle$-orbit of $x$ is nilpotent of class $l$.
 \end{lemma}

\begin{proof}
 Consider the sets
 $$
 Z_s=\{y\in G\mid \langle y^{\langle\f\rangle}\rangle \text{ is nilpotent of class }  s\text{ for every }\f\in A\}.
 $$
 Every set $Z_s$ is closed, because it is equal to the intersection of the closed sets
 $$
 \{y\in G\mid [y^{\f ^{k_0}},y^{\f ^{k_1}},\dots ,y^{\f ^{k_s}}]=1\}
 $$
 over all $\f\in A$ and all tuples $(k_0,k_1,\dots ,k_s)$, where $k_i=0,1,\dots,|\f|-1$. These sets in turn are closed because every automorphism $\f\in A$ is  continuous, as is the operation of taking commutators. Since $A$ is finite, every element $y\in G$ belongs to one of the sets $Z_i$: namely, $y\in Z_m$, where $m$ is the maximum of the nilpotency classes of the subgroups  $\langle y^{\langle\f\rangle}\rangle$ over $\f\in A$. Thus,
$$
G=\bigcup_{s=1}^\infty Z_s.
$$
 Hence by the Baire category theorem \cite[Theorem~34]{kel}, one of these subsets $Z_l$ contains an open subset. This means that there is a coset $y_0K$ of an open subgroup $K$ such that $\langle y^{\langle\f\rangle}\rangle$ is nilpotent of class $l$
for all $\f\in A$ and all $y\in y_0K$. Since $K$ has finite index in $G$, there is a positive integer $N_3$ such that $K$ contains the Sylow subgroups $G_p$ for all $p>N_3$. The coset representative $y_0$ can be chosen in the direct product $\prod_{p\leq N_3}G_p$. Therefore  for $p>N_3$ the projection onto $G_p$ of the nilpotency law of class $l$ for elements of  $\langle y^{\langle\f\rangle}\rangle$ when  $y\in x_0K$  means that $\langle x^{\langle\f\rangle}\rangle$ is nilpotent of class $l$  for all $x\in G_p$ and all $\f\in A$.
\end{proof}

In view of Lemmas~\ref{l-uni-e}, \ref{l-uni-c}, and \ref{l-uni-z},  the following proposition will complete the proof of Theorem~\ref{t-2} in the case where $A$ has prime order.

\begin{proposition}\label{p-unif}
There are functions $N_4(d,q,k,l)$ and $f(d,q,k,l)$ of four positive integer variables $d,q,k,l$ with the following property. Let $p $ be  a prime, and suppose that $P$ is a $d$-generated pro-$p $ group admitting an automorphism $\varphi$ of prime order $q\ne p $ such that $C_P(\varphi )$ is nilpotent of class $k$ and every subgroup  $\langle x^{\langle\f\rangle}\rangle$ generated by the $\langle\f\rangle$-orbit of an element $x\in P$ is nilpotent of class $l$. If $p >N_4(d,q,k,l)$, then  the group $P$ is nilpotent of class at most  $n(d,q,k,l)$.
\end{proposition}

\begin{proof}
It is sufficient to obtain a bound for the nilpotency class in terms of $d$, $q$, $k$,~$l$ for every finite quotient $T$ of $P$ by a $\varphi$-invariant open normal subgroup. Let $L(T)$ be the associated Lie ring constructed with respect to the lower central series of $T$. This Lie ring is nilpotent of exactly the same nilpotency class as $T$ (see, for example, \cite[Theorem~6.9]{khu2}). Therefore it is sufficient to obtain a required bound for the nilpotency class of $L(T)$. Note that all the hypotheses of Proposition~\ref{p-unif} hold for $T$ and its automorphism induced by $\f$, which is denoted by the same letter.

By Proposition~\ref{p-adn-u}, $L(T)$ is generated by $d$ elements such that all commutators $\ka$ in these elements are ad-nilpotent of uniformly bounded index:  $[L(G),{}_{f}\ka ]=0$ for some $f=f(q,k,l)$ depending only on $q,k,l$.
The fixed-point subring $C_{L(T)}(\f)$ is nilpotent of class $k$ because $C_P(\varphi )$ is nilpotent of class $k$.
Thus,  the Lie ring $L(T)$ and its group of automorphisms $\langle \varphi \rangle$  satisfy the hypotheses of Proposition~\ref{p-danilo}. By this proposition, there exist positive integers $e$ and $r$ depending only on $d$, $q$, $k$,   $l$ such that $e\gamma _r(L(T)) = 0$. The additive group of $L(T)$ is a $p $-group. Therefore,  if $p >e$, then $e$ is invertible in the ground ring, so that then $\gamma _r(L) = 0$. It remains to put $N_4(d,q,n,c)=e$ and $n(d,q,n,c)=r-1$.

We thus proved that for $p >N_4(d,q,n,c)$ every finite quotient of $P$ by a $\varphi$-invariant normal open subgroup is nilpotent of class at most $n(d,q,n,c)$. Therefore $P$ is nilpotent of class at most $n(d,q,n,c)$ if $p >N_4(d,q,n,c)$.
\end{proof}

We finally combine all the results in the proof of the main Theorem~\ref{t-2}, which we reproduce here for convenience of the reader.

\begin{theorem}\label{t-2p}
Suppose that a profinite group $G$ admits a finite group $A$ of coprime automorphisms such that for any $\f\in A$  every subgroup $\langle x^{\langle\f\rangle}\rangle$ generated by the $\langle\f\rangle$-orbit of an element $x\in G$ is nilpotent.
\begin{itemize}
  \item[{\rm (a)}] Then the commutator subgroup $[G,A]$ is pronilpotent.
  \item[{\rm (b)}] If in addition $A$ is soluble and the fixed-point subgroup $C_G(A)$ is locally nilpotent, then $G$ is  locally nilpotent.
\end{itemize}
\end{theorem}

\begin{proof} (a)
The commutator subgroup $[G,A]$ is pronilpotent by Theorem~\ref{t-fin}(a) applied to the finite quotients of $G$ by $A$-invariant normal open subgroups.

(b) First we note that it is sufficient to prove the theorem for the case where $A$ has prime order.
  Indeed,  since $A$ is soluble, there is a normal subgroup $A_1$ of prime index $q$. Let $\f\in A\setminus A_1$. Then  $C_G(A_1)$ admits the cyclic group of automorphisms of order $q$ induced by $\langle\f\rangle$, which satisfies on $C_G(A_1)$ the same hypotheses as $A$. Assuming that the theorem is proved in the case where $A$ has prime order, we obtain that $C_G(A_1)$ is nilpotent. Now the group $G$ admits the group of automorphisms $A_1$ of smaller order satisfying the same hypotheses, and therefore the proof concludes by the obvious induction on  $|A|$.   Thus, from now on we assume that $G=\langle\f\rangle$ is cyclic of prime order $q$.

The group $G$ is pronilpotent by  Theorem~\ref{t-fin}(b) applied to the finite quotients of $G$ by $A$-invariant normal open subgroups. Any finite subset $S\subseteq G$ is contained in the $\varphi$-invariant finite set $S^{\langle\varphi \rangle}=\{s^{\varphi ^k}\mid s\in S,\; k=0,1,\dots ,|\f |-1\}$. Therefore we can assume that $G$ is (topologically) finitely generated, say,
 by $d$ elements, and then we need to prove that $G$ is nilpotent. Since $G$ is pronilpotent, we only need to prove that all Sylow $p$-subgroups of $G$ are nilpotent of class bounded by some number independent of~$p$.

Let $N_2$ and $k$ be the numbers given by Lemma~\ref{l-uni-c}, and $N_3$ and $l$ the numbers given by Lemma~\ref{l-uni-z}. Then for every $p>\max\{N_2,N_3\}$, in the Sylow $p$-subgroup $P$, we see that $C_P(\f)$ is nilpotent of class $k$, and every subgroup  $\langle x^{\langle\f\rangle}\rangle$ generated by the $\langle\f\rangle$-orbit of an element $x\in P$ is nilpotent of class $l$. If in addition $p>N_4$ for the number $N_4$ given by Proposition~\ref{p-unif}, then by this proposition $P$ is nilpotent of class  $f(d,q,n,c)$. Since every Sylow $p$-subgroup is nilpotent of some class $c_i$ by Theorem~\ref{t-pro-p}, we obtain a required uniform bound $$\max\{f(d,q,n,c),\, c_i \mid i\leq \max\{N_2,N_3, N_4\}\}$$  independent of $p$  for the nilpotency classes of all Sylow $p$-subgroups. Hence $G$ is nilpotent.
\end{proof}


\begin{thebibliography}{99}

\bibitem{aks} C.~Acciarri, E.~I.~Khukhro, and  P.~Shumyatsky,  Profinite groups with an automorphism whose fixed points are right Engel, \textit{Proc. Amer. Math. Soc.} \textbf{147}, no.~9 (2019), 3691--3703.

\bibitem{acc-shu} C.~Acciarri and  P.~Shumyatsky,
Profinite groups and the fixed points of coprime automorphisms, \emph{J.~Algebra} \textbf{452} (2016), 188--195.

\bibitem{ba-za} Yu. A. Bahturin and M. V. Zaicev, Identities of graded algebras, \textit{J. Algebra} \textbf{205} (1998), 1--12.

\bibitem{bu-me} R. G. Burns and Yu. Medvedev,
A note on Engel groups and local nilpotence, \emph{J.~Austral. Math. Soc. (A)} \textbf{64}, no.~1 (1998), 92--100.

\bibitem{gor}  D. Gorenstein,
\textit{Finite groups},  2nd ed., Chelsea, 1980.

\bibitem{gru}
K. W. Gruenberg, The Engel structure of linear groups, \textit{J.~Algebra} \textbf{3} (1966), 291--303.

\bibitem{gur-kan} R. M. Guralnick and W. M. Kantor,
Probabilistic generation of finite simple groups, \emph{J.~Algebra} \textbf{234}, no.~2 (2000), 743--792.

\bibitem{hal} P. Hall, Some sufficient conditions for a
group to be nilpotent,
\textit{Illinois J. Math.} {\bf 2} (1958), 787--801.

\bibitem{hig}
G. Higman, Groups and Lie rings having automorphisms without non-trivial
fixed points, \textit{J. London Math. Soc.} \textbf{32} (1957), 321--334.


\bibitem{kel} J. L. Kelley, {\it General topology}, Grad. Texts in Math., vol.~27, Springer, New York,  1975.

    \bibitem{khu2} E. I. Khukhro, \textit{$p$-Automorphisms of finite $p$-groups},  London Math. Soc. Lecture Note Ser., vol.~246, Cambridge Univ. Press, 1998.

      \bibitem{kh-sh99}    E. I. Khukhro and  P. Shumyatsky,  Bounding the exponent of a finite group with  automorphisms, \textit{J.~Algebra} \textbf{212}, no.~1 (1999),  363--374.

      \bibitem{kh-sh22}   E.~I.~Khukhro and  P.~Shumyatsky,
       On profinite groups with automorphisms whose fixed points have countable  Engel sinks, \textit{Israel J.~Math.}, \textbf{247}, No. 1 (2022), 303--330.

    \bibitem{laz} M. Lazard, Groupes analytiques $p$-adiques, \textit{Publ. Math. Inst. Hautes \'Etudes
Sci.} \textbf{26} (1965), 389--603.

\bibitem{lin} V. Linchenko, Identities of Lie algebras with actions of Hopf algebras, \textit{Comm.
Algebra} \textbf{25} (1997), 3179--3187.

\bibitem{lu-ma}
A. Lubotzky and A. Mann, Powerful $p$-groups. II. $p$-Adic analytic groups, \textit{J. Algebra} \textbf{105} (1987), 506--515.

\bibitem{rib-zal} L. Ribes and P. Zalesskii, \emph{Profinite groups},
Springer, Berlin, 2010.

\bibitem{rob} D. J. S. Robinson, A course in the theory of groups, Springer, New York, 1996.

\bibitem{rod-shu}
S. Rodrigues and P. Shumyatsky,
Exponent of a finite group admitting a coprime automorphism, \emph{J.~Pure Appl. Algebra} \textbf{224}, no.~9 (2020), Article ID 106370, 8 p.

\bibitem{ser} J.-P.~Serre, \emph{Galois Cohomology}, Springer, Berlin--Heidelberg, 1997.

\bibitem{sha} A. Shalev, Polynomial identities in graded group rings, restricted Lie algebras and $p$-adic analytic groups, \textit{Trans. Amer. Math. Soc.} \textbf{337}, no.~1 (1993), 451--462.

\bibitem{shu00} P. Shumyatsky, Applications of Lie ring methods to group theory, in: \textit{Nonassociative
algebra and its applications} (Eds. R. Costa et al.), Marcel Dekker,
New York, 2000, 373--395.

\bibitem{shu00a} P. Shumyatsky, On pro-$p$
 groups admitting a fixed-point-free automorphism, {\it  J.~Algebra} \textbf{228}, no.~1 (2000), 357--366.

\bibitem{danilo} P. Shumyatsky and D. San\c{c}\~{a}o da Silveira,
On finite groups with automorphisms whose fixed points are Engel,
\textit{Arch. Math.} \textbf{106}, no.~3 (2016), 209--218.

\bibitem{tho}
J. G. Thompson, Finite groups with fixed-point-free automorphisms of prime
order, \textit{Proc. Nat. Acad. Sci. USA} \textbf{45} (1959), 578--581.

\bibitem{wil} J. S. Wilson, \emph{Profinite groups}, Clarendon Press, Oxford, 1998.

\bibitem{wi-ze} J. S. Wilson and E. I. Zelmanov, Identities for Lie algebras of pro-$p$ groups, \emph{J. Pure Appl. Algebra} \textbf{81}, no.~1 (1992), 103--109.

\bibitem{ze92} E. Zelmanov, Nil rings and periodic groups, \textit{Korean Math. Soc. Lecture
Notes in Math.}, Seoul, 1992.

\bibitem{ze95}    E. Zelmanov, Lie methods in the theory of nilpotent groups, in: \textit{Groups'\,93 Galaway/St Andrews}, Cambridge Univ. Press, Cambridge, 1995, 567--585.

\bibitem{ze17} E. Zelmanov,
Lie algebras and torsion groups with identity, \textit{J.~Comb. Algebra} \textbf{1}, no.~3 (2017), 289--340.
 \end{thebibliography}
\end{document}